\documentclass[12pt]{article}

\usepackage{color}

\newfam\msbfam
\font\tenmsb=msbm10  scaled \magstep1 \textfont\msbfam=\tenmsb
\font\sevenmsb=msbm7  scaled \magstep1 \scriptfont\msbfam=\sevenmsb
\font\fivemsb=msbm5  scaled \magstep1 \scriptscriptfont\msbfam=\fivemsb
\def\Bbb{\fam\msbfam \tenmsb}

\def\RR{{\Bbb R}}
\def\CC{{\Bbb C}}
\def\QQ{{\Bbb Q}}
\def\NN{{\Bbb N}}
\def\ZZ{{\Bbb Z}}
\def\II{{\Bbb I}}
\def\TT{{\Bbb T}}

\newtheorem{theorem}{Theorem}

\newtheorem{definition}{Definition}
\newtheorem{example}[definition]{EXAMPLE}

\def\ss{\subseteq}
\def\Aut{\hbox{Aut}}

\def\O{\Omega}
\def\zbar{\overline{z}}
\def\Re{\hbox{Re}\,}
\def\Im{\hbox{Im}\,}
\def\dbar{\overline{\partial}}
\def\zbar{\overline{z}}
\def\0{{\bf 0}}
\def\H{{\cal H}}
\def\F{{\cal F}}
\def\ftilde{\widetilde{f}}

 \def\HollowBox #1#2{{\dimen0=#1 \advance\dimen0 by -#2
       \dimen1=#1 \advance\dimen1 by #2
        \vrule height #1 depth #2 width #2
        \vrule height 0pt depth #2 width #1
        \llap{\vrule height #1 depth -\dimen0 width \dimen1}%
       \hskip -#2
       \vrule height #1 depth #2 width #2}}
 \def\BoxOpTwo{\mathord{\HollowBox{6pt}{.4pt}}\;}

\def\endpf{\hfill $\BoxOpTwo$}

\begin{document}

\begin{center}
{\LARGE \bf Automorphism Groups of Domains}
\medskip \\
{\LARGE \bf that Depend on Fewer Than the}
\medskip \\
{\LARGE \bf Maximal Number of Parameters}
\end{center}
\vspace*{.2in}

\begin{center}
by Jisoo Byun and Steven G. Krantz\footnote{The second author thanks the American Institute of Mathematics
for its hospitality and support during a portion of this work.  Both authors thank the
Banach Center in Warsaw, Poland for its hospitality during a recent conference during
which this problem was discussed.  The second author was supported in part by a grant
from the Dean of the Graduate School at Washington University and a grant from the National
Science Foundation.}
\end{center}

\section{Introduction}

Let $\Omega \subseteq \CC^n$ be a domain, that is a connected,
open set.  Let $\Aut(\Omega)$ denote the collection of biholomorphic
self-mappings of $\Omega$ (see [ISK1] or [GRK] for a survey of this topic).
This set forms a group when equipped with the binary operation of
composition of mappings.  In case $\Omega$ is bounded, then the group
is in fact a real (never a complex) Lie group.   We call this group
the {\it automorphism group} of $\Omega$.

It is naturally a matter of considerable interest to describe the automorphism
group for a given domain $\Omega$.  In the best of all possible circumstances,
we would like to give an explicit description of this group.  As an instance,
in case
$$
\O = B \equiv \{(z_1, \cdots, z_n) \in \CC^n: |z_1|^2 +\cdots +  |z_n|^2 < 1\} \, ,
$$
then the automorphism group of $\Omega$ is generated by {\bf (i)} the
unitary rotations and {\bf (ii)}  the M\"{o}bius transformations 
$$y
(z_1, z_2, \dots, z_n) \longmapsto \left ( \frac{z_1 - a}{1 - \overline{a}z_1} ,
   \frac{\sqrt{1 - |a|^2} z_2}{1 - \overline{a}z_1}, \dots,
       \frac{\sqrt{1 - |a|^2} z_n}{1 - \overline{a}z_1} \right )
$$
for $a \in \CC$, $|a| < 1$.
It is worth noting (as this is part of the theme of the present paper)
that any description of the automorphisms of $B$ {\it must} involve all
$n$ variables.  Such a statement is already true for the unitary group alone.
It is {\it not true} that the automorphism group depends only on $n-1$ variables or
$n-1$ parameters.\footnote{Clearly these statements can be formulated in terms
of the dimension of the automorphism group.  We leave this task as an exercise
for the reader.  For the purposes here, the formulation in terms of the ambient
complex variables is more convenient and more accurate.}

There are other important domains---ones that are currently a
focus of considerable study---whose automorphism groups are
much simpler. For the Kohn-Nirenberg domain (see [JIS],
[KON]), the automorphisms consist of {\it rotations in one
variable only}. The purpose of the present paper is to produce
a geometric criterion which will guarantee that the
automorphism group of a given domain $\O \ss \CC^n$ will
depend on fewer than the full number of variables in the
ambient space. This result will simplify the example in [JIS],
and will also provide further examples for the future. We
indicate some of these latter examples at the end of the
present paper.\footnote{In the paper [GIK], it is shown that
if $\Omega$ is a hyperbolic domain and the dimension of the
automorphism group exceeds $n^2 + 2$ then the domain must be a
ball. The result of the present paper is in a philosophically
similar spirit, but the restrictions on the dimension are more
severe.}

We are happy to thank Peter Pflug for a helpful conversation.

\section{Basic Ideas}

Fix a pseudoconvex domain $\Omega \ss \CC^2$ with smooth
boundary. Because of the discussion in the first section, we
may suppose that $\O$ is not (biholomorphically) the ball. In
case $\Omega$ is strongly pseudoconvex, we may then conclude
by a theorem of Bun Wong and Rosay (see [KRA1]) that
$\Aut(\Omega)$ is compact (in fact we shall make this a
standing hypothesis in the discussion that follows).

\newfam\msbfam
\font\tenmsb=msbm10  scaled \magstep1 \textfont\msbfam=\tenmsb
\font\sevenmsb=msbm7  scaled \magstep1 \scriptfont\msbfam=\sevenmsb
\font\fivemsb=msbm5  scaled \magstep1 \scriptscriptfont\msbfam=\fivemsb
\def\Bbb{\fam\msbfam \tenmsb}

\def\RR{{\Bbb R}}
\def\CC{{\Bbb C}}
\def\QQ{{\Bbb Q}}
\def\NN{{\Bbb N}}
\def\ZZ{{\Bbb Z}}
\def\II{{\Bbb I}}
\def\TT{{\Bbb T}}

Let $(z,w)$ be the coordinates in $\CC^2$. We take it that the
origin $\0 = (0,0)$ lies in $\partial \Omega$. Now suppose that,
near $\0$, $\O$ is defined by the inequality
$$
\rho(z,w) \equiv \Re w + \varphi(z, \overline{z}) + \psi(z, \zbar, \Im w) \cdot (\Im w)^2 < 0 \, . \eqno (*)
$$
In what follows we shall always use the letter $\rho$ to denote a defining
function for $\Omega$.
Applying a unitary rotation if necessary, we may arrange that
the complex tangent space $\H_\0$ at \0 is equal to $\{(z,0): z \in \CC\}$.

Let $\Aut_\0(\O)$ denote the automorphisms of $\O$ that fix $\0$.
We shall assume below that ${\cal H}_0 \cap  \Omega$ is an open subset of ${\cal H}_0$.
If $f$ is an automorphism of $\O$ then set $\widetilde{f}(\zeta) = f_1(\zeta, 0)$
for $(\zeta, 0) \in \H_\0 \cap \O$ and $f =(f_1, f_2)$.  
Define ${\cal F}_{\widetilde{f}}(z,w) = (\widetilde{f}(z), w)$.

Now we shall list the standing hypotheses that will be in place for the
remainder of this paper:

\begin{center}
\bf Standing Hypotheses
\end{center}

\begin{enumerate}
\item[{\bf 1.}] The set $\H_\0 \cap \O$ is  an open subset of $\H_\0$.
\item[{\bf 2.}]  The function $\varphi$ in the definition of $\rho$ has
no harmonic terms.
\item[{\bf 3.}]	 If $f \in \Aut_\0(\O)$ then ${\cal F}_{\widetilde{f}} \in \Aut_\0(\O)$.
\item[{\bf 4.}]  The automorphism group $\Aut_\0(\O)$ is compact.
\item[{\bf 5.}]  The domain $\O$ is complete hyperbolic.
\item[{\bf 6.}]  Any automorphism of $\O$ continues analytically to a neighborhood of $\0 \in \partial \O$.
\end{enumerate}
\vspace*{.15in}

Standing Hypothesis {\bf 2} perhaps merits some discussion.  In case $\varphi$ {\it has}
harmonic terms, we may write
$$
\varphi(z, \overline{z}) = \sum_{k=2}^\infty a_k z^k + \sum_{k=2}^\infty \overline{a}_k \overline{z}^k +
           \widetilde{\varphi}(z, \overline{z}) \, ,
$$
where $\widetilde{\varphi}$ has no harmonic term.  Let $\mu(z) = \sum_k a_k z^k$ be holomorphic.
Then the holomorphic coordinate change
$$
(\widetilde{z}, \widetilde{w}) = (z, w + 2 \sum_k a_k z^k)
$$
defines a new local defining function
$$
\Re \widetilde{w} + \widetilde{\varphi}(z, \overline{z}) + \psi (z, \overline{z}, \Im w) \cdot (\Im w)^2 < 0 \, .
$$
Note that the lead term here has no harmonic term.

Now our main result is this:

\begin{theorem} \sl
Let $\Omega$ be a domain as described above.  Then any
automorphism $f(z) = (f^1(z), f^2(z))$ fixing the origin \0 of $\Omega$ must have
the form
$$
f(z) = (\varphi(z_1), z_2) \, .
$$
In other words, any automorphism of $\Omega$ fixing the origin will depend {\it only}
on the first variable (and not on the second).
\end{theorem}

If $f \in \Aut(\O)$ then certainly $\rho \circ f$ is also a local
defining function for $\O$ near $\0$.  Therefore
$$
\rho \circ f (z,0) = \mu \cdot \rho(z,0)  \eqno (**)
$$
for some positive function $\mu$ near \0.  As a result, we may conclude (see Standing
Hypothesis {\bf 2}) that the quadratic part of $\mu \cdot \rho(z,0)$ has no
harmonic terms.  This observation also applies to the lefthand side of $(**)$,
so we see that
$$
\Re f_2(z,0) + \varphi(f_1(z,0), \overline{f_1(z,0)}) +
     \psi \bigl ( f_1(z,0), \overline{f_1(z,0)}, \Im f_2(z,0) \bigr ) \cdot (\Im f_2(z,0) )^2
$$
has no harmonic terms.
\medskip \\

\noindent {\bf CLAIM:}  We assert that $f_2(z,0) \equiv 0$.
\smallskip \\

\noindent {\bf Proof of the Claim:}  If not, then
$$
f_2(z,0) = \sum_{k \geq 2} a_k z^k \, .
$$
Let $a_{k_0}$ be the nonzero coefficient with least index.

Since $\varphi$ has no harmonic terms and
$$
f_1(z,0) = b_1 z + \sum_{k \geq 2} b_k z^k
$$
with $b_1 \ne 0$, we see that
$$
\varphi(f_1, \overline{f}_1) = \sum_{k,\ell} b_k \overline{b}_{\ell} z^k \zbar^\ell \, .
$$
Certainly $b_0 = 0$, hence $\varphi(f_1, \overline{f}_1)$ has no harmonic terms.

Now
$$
\psi \bigl (f_1(z,0), \overline{f_1(z,0)}, \Im f_2(z,0) \bigr ) \cdot (\Im f_2(z,0) )^2
$$
has the term with index $2k_0$ as the first nonvanishing term.  Hence the first
harmonic term
$$
\Re a_{k_0} z^k \equiv 0 \, .
$$
But this implies that $a_{k_0} \equiv 0$.  And that is a contradiction.
\endpf
\smallskip \\

Now, as a consequence of the claim, we certainly know that
$$
f(\H_\0 \cap \O) \ss \H_\0 \, .
$$
As a result, $f \bigr |_{\H_\0 \cap  \O}$ is an automorphism of
$\H_\0 \cap  \O$ that fixes \0.	In our earlier notation,
$\ftilde \in \Aut(\H_\0 \cap \O)$.

Referring to Standing Hypothesis {\bf 3}, we now consider
$$
\F_{\widetilde{f^{-1}}} \circ f \in \Aut_\0(\Omega) \, .
$$
We see that
\begin{eqnarray*}
\F_{\widetilde{f^{-1}}} \circ f(z,0) & = & \F_{\widetilde{f^{-1}}} (\ftilde(z), 0) \\
                                     & = & (z, 0) \, .
\end{eqnarray*}
Writing $\Phi \equiv \F_{\widetilde{f^{-1}}} \circ f$, we may say
that $\Phi(z,0) = (z,0)$.

Now we have
$$
d\Phi(z,0) = \left ( \begin{array}{cc}
			1 & a(z) \\
			0 & b(z) \\
		   \end{array}
	     \right )
$$  
for some $a(z), b(z)$ that are holomorphic on $\H_\0$. Of
course $d\Phi$ takes the real tangent space at \0 to the real
tangent space at \0 and the complex tangent space at \0 to the
complex tangent space at \0. Since $\Aut_\0(\O)$ is compact, we
conclude that $\sup |b(z)| = 1$ and $b(0) = 1$ hence (by the
maximum principle) $b \equiv 1$.  Finally, since $\O$ is complete
hyperbolic, $a(z) \equiv 0$.  We conclude that
$$
d\Phi(z, 0 ) = \left ( \begin{array}{cc}
                          1 & 0 \\
			  0 & 1 \\
		       \end{array}
	       \right ) \ ,
$$
the identity matrix.   
%
By the Cartan Uniqueness Theorem, we conclude that $\Phi(z) = z$ for all $z \in \Omega$.

Therefore, we finally arrived that 
$$
f(z,w) = \bigl ( \ftilde(z), w) \, .
$$
In short, $f$ depends on fewer than the maximal number of parameters.

\subsection{Levi Flat Domains}

Now we present a variant of our main result.

\begin{center}
\bf Standing Hypotheses
\end{center}

\begin{enumerate}
\item[{\bf 1.}] The local defining function for $\Omega $ at the origin is $\Re w =0$.
\item[{\bf 2.}] Any automorphism of $\O$ continues analytically to a neighborhood of $\0 \in \partial \O$.
\item[{\bf 3.}]	The automorphism group $\Aut_\0(\O)$ is compact.
\item[{\bf 4.}] The domain $\O$ is complete hyperbolic.
\item[{\bf 5.}] If $f \in \Aut_\0(\O)$ then ${\cal F}_{\widetilde{f}} \in \Aut_\0(\O)$.
\end{enumerate}
\vspace*{.15in}

Fix $f \in \Aut_\0(\O)$.  By the same argument as in the Main Theorem, we have that 
$$ \rho \circ f (z,0) = \mu \rho(z,0), $$ where $\rho = \Re w $. This implies that $\Re f_2 \equiv 0$ where
$f=(f_1 , f_2)$. We obtain that $f$ preseve that $H_0\cap \Omega$. By Standing Hypotheses 5, we consider $\Phi = F_{\tilde f^{-1}} \circ f$. Then $\Phi (z,0) = (z,0)$.

By Standing Hypotheses 3, we obtain information on one jet of $f$: 
$$
d\Phi(z, 0 ) = \left ( \begin{array}{cc}
                          1 & 0 \\
			  0 & 1 \\
		       \end{array}
	       \right ) \ .
$$ 
By the power series expansion of $\Phi$, 
\begin{eqnarray*}
\Phi_1 (z,w) &=& z+ w^2  \sum a_{jk} z^j w^k \\
\Phi_2 (z,w) &=& w + w^2\sum b_{jk} z^j w^k .
\end{eqnarray*}

By Standing Hypotheses 1, there is a positive $ \delta$ such that $(z, it)$ is in the boundary of $\O$, for all 
$|z| < \delta$ and  real number $|t| < \delta$. So, we get $\Phi(z,i t)$ lies in the boundary of $\O$. This implies that 
$$\Re \left ( it + (it)^2 \sum b_{jk} z^j (it)^k \right ) = 0 .$$
We consider the left hand side as a power series with respect to $t$.  For all $j, k \ge 0$, 
$$ \Re ( -b_{jk} z^j i^k ) =0. $$

For each $j, k$, we can choose two different complex numbers $\alpha, \beta$ such that $\alpha^j$ is pure 
imaginary  and $\beta^j$ is real. Note that $j$ and $k$ are fixed. 
We can assume that $i^k$ is a real number $\pm1$. We get that 
$$ \Re ( -b_{jk} z^j) =0. $$ 
The above equation holds for all sufficient small complex numbers. We can insert $\alpha, \beta$ into it.
This implies that $b_{jk}$ is a real and pure imaginary number. Hence, all $b_{jk}$ are zero. 
We arrive at the conclusion that 

\begin{eqnarray*}
\Phi_1 (z,w) &=& z+ w^2  \sum a_{jk} z^j w^k \\
\Phi_2 (z,w) &=& w .
\end{eqnarray*}

We want to prove that $\Phi_1 (0,w) = w^2 \sum a_{0k} w^k$ is indentically zero.
Expecting a contradiction, we assume that $a_{0k}$ is the first nonzero term.
We consider the $N$-times composition of $\Phi$ restricted to $(0,w)$. By a calculation,
$$ \Phi^N(0,w) = ( N a_{0k} w^k + \mathrm{(higher order terms)}, w). $$ 
Since $\Phi^N$ is a precompact family, we arrive at a contradiction. 

Therefore, $\Phi(0,w) = (0,w)$. Take a derivative at $(0,w)$. Then 
$$
d\Phi(0, w ) = \left ( \begin{array}{cc}
                          a(w) & b(w) \\
			  0 & 1 \\
		       \end{array}
	       \right ) \ .
$$ 
By the same technique as in the last section, we have the identity matrix at $(0,w)$. We can apply Cartan's Uniqueness Theorem.
Finally we get that $\Phi$ is the identity map of $\O$. This means that $f(z,w) = (\tilde f (z), w )$.

\section{Examples}

\begin{example} \rm
First let us look at the Kohn-Nirenberg domain [KON], which is given
(in the complex variables $z, w$) by
$$
\Re w + |zw|^2 + |z|^8 + \frac{15}{7} |z|^2 \Re z^6 < 0 \, .
$$
It is a simple matter to verify that this domain satisfies the
hypothesis about $S$ in our Theorem.  Also the domain is of finite
type, so every automorphism extends smoothly to the boundary.
We may conclude immediately that any automorphism fixing the origin depends on just
one variable.  Then simple calculations show (see [JIS])
that the only automorphisms are rotations in the $z$ variable.
\endpf
\end{example}

\newfam\msbfam
\font\tenmsb=msbm10  scaled \magstep1 \textfont\msbfam=\tenmsb
\font\sevenmsb=msbm7  scaled \magstep1 \scriptfont\msbfam=\sevenmsb
\font\fivemsb=msbm5  scaled \magstep1 \scriptscriptfont\msbfam=\fivemsb
\def\Bbb{\fam\msbfam \tenmsb}

\def\RR{{\Bbb R}}
\def\CC{{\Bbb C}}
\def\QQ{{\Bbb Q}}
\def\NN{{\Bbb N}}
\def\ZZ{{\Bbb Z}}
\def\II{{\Bbb I}}
\def\TT{{\Bbb T}}

\begin{example} \rm
Let $\varphi$ be a $C^\infty$ function on $\RR$,  even, nonnegative, supported
in the interval $[-1/10, 1/10]$, constantly equal to $1$ on $[-1/100, 1/100]$.
Define
\begin{eqnarray*}
\Omega & = & \{(z_1, z_2) \in \CC^2: \bigl (1 - \varphi(|z_2|^2) \bigr )\bigl (1 - \varphi(|1 - z_1|^2) \bigr ) \cdot [-1 + |z|^2] \\
    && \quad + \varphi(|z_2|^2)\varphi(|1 - z_1|^2) \cdot [- 1 + \epsilon + \Re z_1] < 0\} \, .
\end{eqnarray*}
Then $\Omega$ is nothing other than the unit ball in $\CC^2$ with a flat bump
centered about the spherical boundary point $(1,0)$.  And notice that the boundary
point $(1-\epsilon, 0)$ has a neighborhood in the boundary that lies in the hyperplane
$\{\Re z_1 = 1-\epsilon\}$.  It is straightforward to see that the only automorphisms of $\Omega$
are rotations in the $z_2$ variable (see [LER]).

Now examine our main theorem. This domain $\Omega$ satisfies
the hypotheses of that theorem with the origin replaced by $(1-\epsilon, 0)$. 
The automorphism group of this particular $\Omega$ may be determined explicitly
(in fact any automorhism of $\Omega$ is just an automorphism of the unit ball---see[LER]),
so most of the Standing Hypotheses are automatic.
Note particularly that the presence of the bump
near $(1-\epsilon,0)$ forces standing hypothesis {\bf 3}\ to hold. So
this example is an illustration of our main result. The
automorphism group only depends on the $z_2$ variable.
\end{example}

\section{Concluding Remarks}

This is the first paper to explore the questions posed here.  There is clearly
a need for a result in all dimensions, and for results that have more flexible
hypotheses.

It would certainly be of interest to have concrete examples to which our results
do (or do not) apply.  The theorem definitely does {\it not} apply to the
complex ellipsoids
$$
E_{m,n} = \{(z_1, z_2): |z_1|^{2m} + |z_2|^{2n} < 1\}
$$
for $m, n \in \NN$.  And of course they should not.  The result also does
not apply to the ball or the Siegel upper half space.  More, they
do not apply to any of the bounded symmetric domains of Cartan (see [GRK] for
a discussion of some of these specialized domains).

What would be ideal is to have a theorem that, given $0 < k < n \in \NN$, characterizes
domains in $\CC^n$ whose automorphisms depend only on $k$ variables.  This will
be the subject of future investigations.
\vfill
\eject

\noindent {\Large \sc References}
\bigskip \\

\begin{enumerate}


\item[{\bf [BEP]}]  E. Bedford and S. Pinchuk, Domains in ${\CC}^2$ with non-compact holomorphic
automorphism group (translated from Russian),
{\it  Math.\ USSR-Sb.} 63(1989), 141--151.\

\item[{\bf [BER]}]  L. Bers, {\it Introduction to Several Complex Variables},
New York Univ. Press, New York, 1964.

\item[{\bf [JIS]}]  J. Byun and H.R. Cho, Explicit Description for the Automorphism Group of
the Kohn-Nirenberg Domain, {\it Math.\ Z.}, to appear.


\item[{\bf [CHES]}]  S.-C. Chen and M.-C. Shaw, {\it Partial Differential
Equations in Several Complex Variables}, American Mathematical Society,
Providence, RI, 2001.

\item[{\bf [JPDA]}]  J. P. D'Angelo, {\it Several Complex Variables and
the Geometry of Real Hypersurfaces}, CRC Press, Boca Raton, 1992.	

\item[{\bf [DIF1]}]  K. Diederich and J. E. Forn\ae ss, Pseudoconvex
domains:  An example with nontrivial Nebenh\"{u}lle, {\it Math. Ann.}
225(1977), 275--292.


\item[{\bf [GIK]}]  J. Gifford, A. Isaev, and S. G. Krantz, On the Dimensions
of the Automorphism Groups of Hyperbolic Reinhardt Domains,
{\it Illinois Jour.\ Math.} 44(2000), 602--618.

\item[{\bf [GRK]}]  R. E. Greene and S. G. Krantz, Biholomorphic self-maps of domains,
{\it Complex Analysis II} (C. Berenstein, ed.), Springer Lecture
Notes, vol. 1276, 1987, 136-207.

\item[{\bf [ISK1]}] A. Isaev and S. G. Krantz, A. Isaev and S.
G. Krantz, Domains with non-compact automorphism group: A
survey, {\it Advances in Math.} 146(1999), 1--38.

\item[{\bf [JAP]}] M. Jarnicki and P. Pflug, {\it Invariant
Distances and Metrics in Complex Analysis}, Walter de Gruyter
\& Co., Berlin, 1993.

\item[{\bf [KOH]}]  J. J. Kohn, Boundary behavior of $\dbar$ on weakly
pseudoconvex manifolds of dimension two, {\em J. Diff. Geom.} 6(1972),
523--542.

\item[{\bf [KON]}]  J. J. Kohn and L. Nirenberg, A pseudo-convex domain
not admitting a holomorphic support function, {\em Math. Ann.} 201(1973), 265-268.

\item[{\bf [KRA1]}]  S. G. Krantz, {\it Function Theory of
Several Complex Variables}, $2^{\rm nd}$ ed., American
Mathematical Society, Providence, RI, 2001.

\item[{\bf [KRA2]}]  S. G. Krantz, Determination of a domain in complex space by its
automorphism group, {\it Complex Variables Theory Appl.} 47(2002), 215--223.

\item[{\bf [LER]}] L. Lempert and L. A. Rubel,
An independence result in several complex variables.
{\it Proc.\ Amer.\ Math.\ Soc.} 113(1991), 1055--1065.
					
\item[{\bf [NIV]}]  I. Niven, {\it Irrational Numbers}, The Mathematical Association
of America in affiliation with John Wiley \& Sons, New York, 1956.

\item[{\bf [RUD]}]  W. Rudin, {\it Function Theory in the Unit Ball
of $\CC^n$}, Grundlehren der Mathematischen Wissenschaften in
Einzeldarstellungen, Springer, Berlin, 1980.

\end{enumerate}

\end{document}